# Première approche de la densité d'un opérateur de Perron Frobenius.
## III – Applications : EDP, EDO, etc.

Guy Cirier[*]

**Abstract**

First approach of invariant densities of a Perron Frobenius operator. Asymptotic behaviours of ODE or PDE, as $\partial a/\partial t = F(a)$, are most interesting. The associed infinitesimal iteration is $f(a,\delta) = a + \delta F(a)$. If $F(a)$ is partially linear, a random distribution can be asymptotic solution. Among applications, are asymptotic profiles of of Lorenz, Navier Stokes or Hamilton's équations.

**Résumé**

À toute EDO ou EDP, de la forme $\partial a/\partial t = F(a)$, on peut associer l'itération infinitésimale $f(a,\delta) = a + \delta F(a)$ et lui appliquer la densité de Perron Frobenius. Si $F(a)$ est partiellement linéaire, alors une distribution aléatoire peut être solution asymptotique de cette équation. Les applications les plus remarquables de l'étude sont les solutions asymptotiques des équations de Lorenz, de Navier Stokes ou d'Hamilton.

**Mots cléfs** : Perron Frobenius et EDP, seuil de récurrence ; équations de Lorenz ou Navier Stokes.
**Key words**: Perron Frobenius and EDP , recurrent doorstep. Lorenz or Navier Stokes equations.

## I – Introduction

Cet article présente quelques applications de la densité invariante [3] [4] et ne peut être lu sans ces préalables. Moins important que les deux premiers, il a été traité un peu plus superficiellement du point de vue des mathématiques. Il intéressera davantage les physiciens ou les mécaniciens. Toutefois, les domaines mathématiques couverts vont de certaines équations intégrales à certaines équations différentielles en passant par des équations avec décalages. Le domaine de cette approche commence là où les techniques déterministes classiques grippent, et où le probabilisme n'a pas encore obtenu complètement droit de cité. Toutefois, certaines questions méritent bien des approfondissements.
Quelques équations (Lorenz, Navier Stokes) montrent la puissance d'investigation de la méthode.

## II- Équations différentielles ordinaires ou aux dérivées partielles

On ne considère ici que les équations de la forme $\partial a / \partial t = F(a)$ où les inconnues forment un vecteur $a$ de $R^d$ et les variables $t = (q, t_k)$ de $R^{+k}$ s'interprètent comme la position $q = (q_1,..,q_i,...q_{k-1})$ supposée ici positive de $R^{+k-1}$ et le temps $t_k \in R^+$. Le cas d'une EDO est celui où $k = 1$ et la seule variable est le temps $t$. La référence sur la question reste Arnold [1] et [2]. Pour étudier une équation différentielle à l'aide de l'opérateur de Perron Frobenius, il faut d'abord traduire l'équation en termes d'itération.

**Définition 1**
*On appelle itération différentielle l'application $f(a,\delta)$ de $R^d$ dans $R^d$, définie par $f(a,\delta) = a + \delta F(a)$ pour $\delta_0 > \delta > 0$ fixé de $R^{+k}$ (ainsi utilisée dans les calculs numériques).*
Notons que les zéros réels de $F$ sont les points fixes de $f(a,\delta)$ et qu'ils sont supposés isolés.

---

[*] guy.cirier@gmail.com

**Hypothèse 0**

*On suppose ici, pour utiliser les résultats de* $[3]$ *et* $[4]$, *que* $f$ *polynomiale itère un compact* $C$ *de* $R^d$ *dans lui-même et que l'on peut appliquer la méthode du col.*

**Difficultés pour une approche probabiliste**

On veut d'abord déterminer les conditions de convergence vers une distribution asymptotique, et préciser ensuite cette distribution. Mais, pour tout $\delta > 0$ fixé, même en admettant que l'on puisse appliquer tous les résultats des articles précédents à $f(a,\delta)$, la distribution intrinsèque sur un domaine invariant peut dépendre de $\delta$. En outre, il faut rappeler que les cycles attractifs $f^{(n)}(a,\delta) = a$ appartenant à $F \circ Fix(f)$ peuvent masquer les éventuelles solutions probabilistes intrinsèques.

On montre alors dans ce qui suit que, lorsque $\delta \to 0$ et que la hessienne est définie, en général le domaine invariant $S(\delta) \to \{0\}$ point fixe ou un cycle. Mais, la dégénérescence de la hessienne peut induire un comportement aléatoire complexe. On étudie alors l'itération $f(a,\delta)$ partiellement linéaire. Rappelons que, dans le cas des EDO, des théorèmes déterministes, du type Poincaré-Bendixon, donnent alors des comportements asymptotiques simples, mais ne fonctionnent plus au-delà de la dimension 2. Si $F$ polynomiale applique un compact dans lui-même, $F$ est lipschitzienne et, en partant d'un point de ce compact, il existe une solution locale. Mais cela ne donne aucune règle sur le comportement quand $t \to \infty$.

**1- Rappels et notations pour une approche probabiliste**

Quitte à introduire des linéarités et des variables supplémentaires ou des relations de compatibilité, on peut ramener beaucoup d'équations à une équation de la forme $\partial a/\partial t = F(a)$, où $F(a)$ est une application de $R^d$ dans $R^d$, avec autant d'équations que d'inconnues où $t \in R^{+k}$ avec $k \leq d$.

L'indice $\ell = 1, 2, ...d$ est réservé aux inconnues $a$ et l'indice $i = 1, 2, ..k$ aux variables $t$. On regroupe par blocs $L_i$ avec $\Sigma_{i \in L_i} \ell_i = d$ les indices $\ell$ des inconnues $a_\ell$ relatifs à une même variable d'indice $i$ dans $\partial a/\partial t$. Ces indices de $a_{\ell,i}$ sont notés $\ell,i$. Soit $\delta_i = t_i / n'$ un partage de $t_i$ en $n'$ parties égales ; il reste constant quand $\ell, i$ parcourt les $L_i$.

$\delta$ petit vecteur de $R^{+k}$ étant fixé, on itère $n \in N^d$ fois $f(a,\delta)$ en faisant tendre $(n',n) \to \infty$ et en partant d'une position initiale $a(0) = (a_\ell(t(0)))$ en $t(0) = (q(0),0)$. Classiquement, la solution $a(t)$ des équations différentielles est obtenue en liant $n' = n$, ce qui rend l'itération $f(a,\delta)$ dépendant fonctionnellement de $n$ et l'on en cherche la limite quand $n \to \infty$.

**Problème de cohérence des notations**

Pour appliquer les méthodes des articles $[3][4]$, la question se pose au sujet de la fréquence des itérations partielles des coordonnées $f_\ell(a,\delta) : (\zeta_\ell) = \lim_{n \to \infty}(n_\ell)/\Sigma_\ell(n_\ell)$. Même s'il semble plus naturel de se placer dans le cas traditionnel $(\zeta_\ell) = 1/d$, il restera utile de considérer la fréquence des itérations des coordonnées pour l'étude de la résonance, envisagée dans le prochain article.

De même, il est possible en théorie de prendre des indices de partage $n' \in N^k$ qui tendent vers $\zeta'_{\ell'} = \lim_{n' \to \infty} n'_{\ell'} / \Sigma_{\ell'} n'_{\ell'}$. Mais, puisque $\delta_i = t_i / n'$, il est plus parlant de faire varier $t$ plutôt que $n'$ pourvu que le ratio $\delta_i$ reste petit. On verra apparaître un problème quand la variable $t \to \infty$ : la convergence asymptotique vers tel ou tel ensemble va dépendre de la direction $\tau$ de $t$ qu'il faut se fixer : $\tau = \lim_{t \to \infty} t / |t|$ avec $|t| = \Sigma_i t_i$. Dès lors, $\delta = \tau |t| / n \to 0$ dans la direction $\tau$ et la distinction entre $n'$ et $n$ disparaît. Dans le cas des EDO ou lorsque la variable $q$ reste bornée, $\tau = 1$.

**Définition 2**

*L'itération différentielle est* $f(a,\delta)$ *dite partiellement linéaire si* $F$ *est partiellement linéaire.*

Comme au chapitre V de [4], on n'étudie ici que ce cas de dégénérescence de la hessienne :
$F$ est partiellement linéaire : $F(a) = (\lambda a + A(a)b + B(a), \lambda'b + C(a)b + D(a))$ où $a = (a,b)$ avec $a \in R^p$ et $b \in R^{d-p}$. La partie linéaire de $F$ est diagonalisée : $A(a)$, $C(a)$ sont au moins de degré 1 et les $B(a), D(a)$ de degré 2 au moins. Pour simplifier, $F$ a deux zéros réels $\alpha = (\alpha, \beta)$ et $0$ isolés.

**Définition 3**
À chaque zéro réel de $F$, on associe une itération asymptotique partiellement linéaire :
$G(a) = (a + \tau A(a)b, b)$ en $0$ et $G_\alpha(u) = (u + \tau A(u+\alpha)v, v)$ au point fixe $\alpha$ avec $a = u + \alpha$.

**Résultat principal : Théorème**
*Si l'itération différentielle $f(a,\delta) = a + \delta F(a)$, de $R^d$ dans $R^d$ pour tout $\delta_0 > \delta > 0$ où $\delta$ de $R^{+k}$, vérifie les hypothèses H0 des articles [3][4] avec pour mémoire:*
*$f(a,\delta)$ est polynomiale et applique un compact $C$ dans lui-même. Les points fixes $0, \alpha \in C$, alors :*
*- si la hessienne de $yF(a)$ n'est pas structurellement dégénérée, outre les points fixes, il n'existe que des orbites cycliques de densité uniforme, de période $T$ telle que $\int_0^T F(a(t+t_0))dt = 0$, $\forall t_0$ ;*
*- si $F$ est partiellement linéaire, outre les points fixes et les orbites cycliques, la solution intrinsèque asymptotiquement invariante lorsque $\delta = |t|\tau/n \to 0$ tend vers celle de $G(a) = (a + \tau A(a)b, b)$ en $0$ et celle de $G_\alpha(a) = (a + \tau A(a+\alpha)b, b)$ en tout point fixe $\alpha$.*
*Selon que $\text{Re}(x\tau A(a)\beta)$ est positif ou négatif au point critique, on prendra la distribution relative à $0$ ou celle de $\alpha$ qui rend contribution de la méthode du col maximum. Les points critiques réalisant $\text{Re}(x\tau A(a)\beta) = 0$ créent une voie de communication entre ces deux ensembles.*

Les étapes de la démonstration de ce théorème suivent. Elles sont guidées par la méthode du col et la géométrie de l'attracteur de Lorenz.

**2- Conditions de convergence**
Il faut étudier les conditions de convergence vers un hypothétique domaine invariant $S(\delta)$. Sous H0, on va déterminer les conditions de convergence pour tout $\delta > 0$ fixé assez petit. Les valeurs propres réelles $\rho$ de la partie linéaire de $f(a,\delta)$ de la forme $\rho \sim 1 + \lambda\delta$ vérifient $|(\partial F(0)/\partial a) - \lambda I| = 0$. Quand $\delta \to 0$ dans la direction $\tau$, les conditions de convergence de [3] vers $S(\delta)$ relatives aux valeurs propres se transposent sur $\lambda\tau$.

**Lemme 1**
*La distribution invariante l'emporte sur le point fixe ou non selon que $\lambda\tau$ est positif ou négatif. Pour une EDO, $\Sigma\lambda > 0$ suffit pour faire tendre la distribution invariante vers $S(\delta)$ s'il n'y a pas de cycles. Dans le cas des EDP, les plans $\lambda\tau = 0$ sont des failles. Si $\tau = 1$, on a la résonance classique.*
Comme dans [3], l'écart résolvant relatif à l'itération $f(a,\delta)$ s'écrit $e_f^n(y,\delta) = y^n - H_n(y,\delta)$ avec $H_n(y,\delta) = \partial^n e^{yf(a,\delta)}/\partial a^n\big|_{a=0}$. Pour tout $\delta > 0$ fixé, $y = ns$ transforme l'ensemble des zéros de $H_n(y,\delta)$ en $S_n(\delta)$. Dans le cas traditionnel où toutes ses composantes d'itération sont égales à $n$, la convergence vers le domaine invariant se traite en posant dans l'écart résolvant $e_f^n(y)$ : $y = ns$ et en mettant $\rho^n$ en facteur : $n^{-n}e_f^n(ns) = s^n - \rho^n H_n(s,\rho)$ doit tendre vers 0 pour tout $s \neq 0$ et tout $n$. Faisons tendre $n$ vers l'infini: $\rho^n = (1+\lambda\delta)^n \sim e^{n\lambda\delta}$ avec $t = \delta = \tau|t|/n$, quand $n \to \infty$ pourvu que $|t|/n$ reste borné. L'expression tend vers l'infini, 0 ou 1 selon que $\lambda\tau = \Sigma_{\ell i}\lambda_{\ell i}\tau_{\ell i}$ est positif, négatif ou nul. Les zéros de $H_n(s,\rho)$ imposent alors ou non leur limite.

Dans le cas des EDP, considérons $\lambda\tau = 0$. En ce cas, il suffit donc de faire varier la direction $\tau$ pour rendre $\lambda\tau$ rationnel et alors passer brutalement d'un point fixe à un ensemble probabilisé et inversement. On a de véritables failles dans le comportement.
Notons que si le point fixe est répulsif, il peut exister des cycles attractifs qui masquent ces problèmes.

**3 – Cycles ou points fixes si la hessienne n'est pas dégénérée.**
Il s'agit maintenant d'étudier les zéros de $H_{n-1}(y,\delta)$. Par la méthode du col $[4]$, l'étude se ramène à celle des points critiques imaginaires de l'application de Plancherel-Rotach avec : $\gamma(a,\delta) = ya + \delta y f(a,\delta) - n\log a$ où $\delta = \tau|t|/n$. Mais on doit considérer tous les points de $F \circ Fix(f)$. En $0 \in Fix(f)$ le domaine $S(\delta)$ est déterminé par les points critiques de $\gamma(a,\delta)$, sous réserve que $\gamma(a,\delta) \to -\infty$ quand $(y,n) \to \infty$ et que sa hessienne soit définie négative.

**Lemme 2**
*Sous les hypothèses de $[4]$, si la hessienne est définie, le point critique $a_c$ est asymptotiquement réel quand $\delta \to 0$.*

- Cas non cyclique ; les équations du point critique $a_c$ de $\gamma(a,\delta)$ sont $a_{c\ell}\partial\gamma a_c,\delta)/\partial a_\ell = y_\ell a_{c\ell} + \delta\partial y F(a_c)/\partial a_\ell - n_\ell = 0$, $\ell = 1,2,..,d$. Si $\delta \to 0$ et $n \to \infty$, $a_{c\ell}\partial\gamma a_c,\delta)/\partial a_\ell \sim y_\ell a_{c\ell} - n_\ell \to 0$ et $a_{c\ell}$ est réel. Il n'y a pas de distribution invariante en dehors des points fixes et des cycles.

- Cas cyclique : $f^{(n)}(a,\delta) = a$. Notons que si $\delta \to 0$ sans dépendre de $n$, les cycles vont disparaître $f^{(n)}(a,\delta) = f^{(n-1)}(a,\delta) + \delta F \circ f^{(n-1)}(a,\delta) = a$, car à la limite $f^{(n)}(a,0) = f^{(n-1)}(a,0) = a$. Toutefois, pour tout $\delta = t/n > 0$, il se peut qu'il y ait des solutions de période $T : \lim_{n\to\infty} f^{(n)}(a,T/n) = a$. Dans le cas d'une EDO, chaque point étant visité un même nombre de fois, la densité de probabilité asymptotique est uniforme. De plus, de la relation de récurrence, il vient $\partial f^{(n)}(a,\delta)/\partial\delta = \left[1 + \delta F \circ f^{(k)}(a,\delta)\right]\partial f^{(n-1)}(a,\delta)/\partial\delta + F \circ f^{(n-1)}(a,\delta)$ et, quand $\delta \to 0$, on a par récurrence pour tout $a$ du cycle: $(1/n)\sum_{k==0}^{k=n-1} F \circ f^{(k)}(a,\delta) \to 0$ qui se traduit asymptotiquement par $\int_0^T F(a(t+t_0))dt = 0, \forall t_0$. Dans le cas d'une EDP, on obtient des résultats similaires.

**Exemples**
**La logistique** $da/dt = \alpha(a - a^2)$. $f(a,\delta) = a + \delta F(a) = a + \delta\alpha(a - a^2)$.
Des 2 points fixes $a = 0$ et $a = 1$, seul $a = 1$ est attractif car le multiplicateur est $1 - \delta\alpha < 1$. La distribution limite est réduite à ce point.

- **k corps de masse $m_i$ dans un champ de forces dérivant d'un potentiel $U(q)$**
Les équations d'Hamilton sont : $dp/dt = -\partial U/\partial q$, $dq/dt = p/m$ de coordonnées $p_i/m_i$ ;
Sous réserve de validité des calculs, la fonction de Plancherel-Rotach de l'itération différentielle associée: $\gamma(p,q) = xp + yq + \delta(-x\partial U/\partial q + yp/m) - n\log p - n\log q$ se décompose en deux fonctions indépendantes : $\gamma(p,q) = \gamma(p) + \gamma(q)$ où $\gamma(p) = (xp + \delta yp/m - n\log p)$ et $\gamma(q) = (yq - \delta x\partial U/\partial q - n\log q)$.
d'où les équations des points critiques $(p,q)$ :
$\partial\gamma/\partial p = x - \delta x\partial^2 U/\partial p\partial q + \delta y/m - n/p = x + \delta y/m - n/p = 0$
$\partial\gamma/\partial q = y - \delta x\partial^2 U/\partial q^2 - n/q = 0$.
Si $\gamma(p,q)$ a une hessienne non dégénérée, il n'existe pas de solution aléatoire. Si le potentiel est newtonien ($U(q)$ n'est pas polynomial), sous réserve de validité des calculs, on a des cycles.

## 4- Itération différentielle partiellement linéaire

Revenons à l'itération différentielle partiellement linéaire $f(a,\delta) = a + \delta F(a)$ avec $\delta = |t|\tau/n$ où $F(a)$ est partiellement linéaire (définition 2). Soient les itérations asymptotiques (définition 3) $G(a) = (a + \tau A(a)b, b)$ en $\mathbf{0}$ et $G(a) = (a + \tau A(a)b, b)$ en $\alpha$ où $a = u + \alpha$.

**Lemme 4**

*L'itération $f(a,\delta) = a + \delta F(a)$ pour $\delta > 0$, où $F(a)$ est partiellement linéaire partant du point $a(0) = (a(0),b(0))$ est asymptotiquement équivalente à l'itération $G(a)$ après transformation d'échelle quand $\delta \to 0$. En $\alpha$, $f(a,\delta)$ est asymptotiquement équivalente à l'itération $G_\alpha(u)$.*

Pour $\delta = |t|\tau/n$, on pose $b_1 = b|t|/n$ et $a_1 = (a,b_1)$ dans $f(a,\delta)$:

$f(a_1,\delta) = (a + \tau A(a)b_1 + \delta(\lambda a + B(a)), b_1 + \delta(\tau\lambda'b_1 + tC(a)b_1) + \delta^2 D(a))$.

Lorsque $\delta \to 0$, c'est à dire quand $n \to \infty$, si reste $b_1$ constant, $f(a_1,\delta) \to (a + \tau A(a)b_1, b_1)$.

Ce phénomène localise les distributions intrinsèques autour de chaque zéro de $F$ et transforme $\mathbf{0}$ en point fixe de $G(a)$. Soit $\alpha = (\alpha,\beta)$ un autre zéro arbitraire de $F$ et distinct de $\mathbf{0}$ : $F(\alpha) = (0,0)$, l'itération translatée (proposition 1 de [3]) de $\alpha$ $f_\alpha(a,\delta) = f(a+\alpha,\delta) - \alpha$ sera encore partiellement linéaire. Et $G_\alpha(a_1) = (a + \tau A(a+\alpha)b_1, b_1)$ sera l'itération asymptotique de $f_\alpha(a,\delta)$ localisée en $\alpha$. On allège l'écriture en posant $b_1 = b$ et $a_1 = a$

**Proposition 1**

*Selon que $\text{Re}(x\tau A(a)\beta)$ est positif ou négatif, on devra prendre la famille de distributions relatives à $\mathbf{0}$ ou à $\alpha$ rendant ainsi la contribution de la méthode du col maximum. Les points critiques réalisant $\text{Re}(x\tau A(a)\beta) = 0$ créent une voie de communication entre ces deux ensembles. Comme $G$ et $G_\alpha$ sont partiellement linéaires, on aura des familles de variétés aléatoires par itération du point critique.*

Pour $\delta = 0$, la distribution asymptotique $G_\alpha(a)$ en $\alpha$ n'est pas la translatée de $G(a)$ en $\mathbf{0}$. En effet (voir proposition 1 de [3]), si $Ta = a + \alpha$, $T^{-1}a = a - \alpha$ et $Tf_\alpha T^{-1} = f$. Si $f(a,\delta) \sim G(a)$ et $f_\alpha(a,\delta) \sim G_\alpha(a)$. On devait donc avoir $TG_\alpha T^{-1} \sim G$. Or, on constate que l'itération calculée $\alpha$ et amenée en $\mathbf{0}$ : $TG_\alpha T^{-1}(a) = (a + \tau A(a)(b-\beta), b) = G(a) - (\tau A(a)\beta, 0)$. Cela s'explique par le fait que $\alpha$ n'est pas point fixe de $G$.

Mais, pour tout $\delta \neq 0$, comme $\alpha$ est point fixe de $f(a,\delta)$, $f(a,\delta)$ et $f_\alpha(a,\delta)$ déterminent un seule et même distribution invariante (proposition 1 de [3]). Pour tout $\delta > 0$ petit, $f(a,\delta) \sim G(a)$ en $\alpha$ ou en $\mathbf{0}$ et $f_\alpha(a,\delta) \sim G_\alpha(a)$ en $\alpha$. Le chapitre IV de [4] permet de lever l'ambiguïté :

Rappelons que, pour tout $\delta > 0$ fixé, on doit considérer la répartition des zéros de $H_{n-1}(y)$, pour obtenir la distribution invariante. Pour $\delta \neq 0$ très petit, on applique la méthode du col au second membre. Il faut alors à sélectionner les points critiques imaginaires de $\gamma(a,\delta) = yf(a,\delta) - n\log a$ qui rendent la partie réelle de $\gamma(a)$ maximum. En prenant $\mathbf{0}$ comme origine et $y = (x, y)$, on a 2 candidats, $yG(a)$ et $yTG_\alpha T^{-1}(a) = yG(a) - x\tau A(a)\beta$, à un infiniment petit près pour représenter $yf(a,\delta)$. Le signe de $\gamma(a,\delta) - \gamma_\alpha(a,\delta) = x\tau A(a)\beta$ donne le maximum.

## 5 - Les équations de Lorenz ( voir [10] par exemple)

Cette équation est particulièrement importante puisque la linéarité partielle de son itération différentielle va générer une solution probabiliste asymptotique. Elle devient ainsi l'exemple idéal pour exposer les résultats de la théorie développée dans les 2 premiers articles.

- *Présentation de l'itération différentielle, de ses points fixes et de leur « répulsivité »*

Ces équations s'écrivent dans nos notations $d\boldsymbol{a}/dt = F(\boldsymbol{a})$, où $\boldsymbol{a} = (a,b,c)$ : ( $da/dt = \sigma(b-a)$ ; $db/dt = \rho a - b - ac$ ; $dc/dt = -\beta c + ab$ ). L'itération différentielle associée ($a + \delta\sigma(b-a)$ ; $b + \delta(\rho a - b - ac)$ ; $c + \delta(-\beta c + ab)$) est quadratique mais présente une linéarité par rapport à $a$. L'itération est traditionnelle et applique un compact dans lui-même pour $\delta_0 > \delta > 0$ (le phénomène se produisant entre une plaque froide et une plaque chaude). Les points fixes sont les zéros de $F(\boldsymbol{a}) = 0$. Lorsque $\rho > 1$ et $\alpha = \sqrt{\beta(\rho-1)}$, il en existe trois : le point $\boldsymbol{0} = (0,0,0)$, et deux autres symétriques par rapport à l'axe $c$ : $\boldsymbol{\alpha}_+ = (\alpha, \alpha, \alpha^2/\beta)$ et $\boldsymbol{\alpha}_- = (-\alpha, -\alpha, \alpha^2/\beta)$.

En $\boldsymbol{0}$, l'équation aux valeurs propres de la partie linéaire est : $(\beta+\lambda)[(\sigma+\lambda)(1+\lambda) - \sigma\rho] = 0$, mais devient en $\boldsymbol{\alpha}_+$ ou en $\boldsymbol{\alpha}_-$ : $\lambda(\beta+\lambda)(1+\sigma+\lambda) - \alpha^2(\lambda + 2\sigma) = 0$. Les $(\beta, \rho, \sigma)$ sont ici supposés tels que ces trois points fixes soient répulsifs, ce qui oblige à étudier les distributions aléatoires autour de chaque point fixe. On oublie ici les cycles attractifs, les résonances éventuelles, pour certaines valeurs des paramètres, etc. Il reste donc beaucoup de choses à préciser.

- *Fonction de Plancherel Rotach au point fixe $\boldsymbol{0}$*

En $\boldsymbol{0}$, posons $\boldsymbol{a}' = \delta\boldsymbol{a}$ qui doit rester borné car la première équation devient $a_1 + \delta\sigma(\delta b - a_1)$ et si $\delta \to 0$ l'itération différentielle tend asymptotiquement vers :

$G(\boldsymbol{a}_1) = (a_1;\ b + \rho a_1 - a_1 c;\ c + a_1 b)$. Soit le vecteur $\boldsymbol{y} = (x,y,z)$. Considérons la fonction de Plancherel Rotach : $\gamma(\boldsymbol{a}_1) = \boldsymbol{y}G(\boldsymbol{a}_1) - nLog a_1 - nLog b - nLog c = \boldsymbol{y}G(\boldsymbol{a}_1) - nLog a_1 bc$

En posant $\boldsymbol{y} = n\boldsymbol{s}$, avec $\boldsymbol{s} = (r,s,t)$ et $\boxed{\varpi = r + s\rho}$, on peut écrire $\boldsymbol{y}G(\boldsymbol{a}_1)/n = a_1\varpi + bs + ct - sa_1 c + tba_1 = L(\boldsymbol{a}_1) + Q(\boldsymbol{a}_1)$

où $L(\boldsymbol{a}_1) = a_1\varpi + bs + ct$ et $Q(\boldsymbol{a}_1) = -sa_1 c + tba_1$

Puisque l'application est quadratique, on va traiter directement les calculs avec des polynômes d'Hermite. Examinons d'abord la hessienne de $G(\boldsymbol{a}_1)/n$ qui est celle de la quadrique $Q(\boldsymbol{a}_1)$ :

$Q = \begin{pmatrix} 0 & t & -s \\ t & 0 & 0 \\ -s & 0 & 0 \end{pmatrix}$ dont la matrice des vecteurs propres est : $T = \dfrac{1}{\mu\sqrt{2}}\begin{pmatrix} 0 & \mu & \mu \\ s\sqrt{2} & -t & t \\ t\sqrt{2} & s & -s \end{pmatrix}$ où $\mu$ est la valeur propre positive de l'équation caractéristique de $Q$ : $\mu(\mu^2 - s^2 - t^2) = 0$. La hessienne est constante et dégénérée et $T$ est orthogonale constante en tout point $\boldsymbol{a}_1$.

- *Changement de base*

L'application orthogonale $\boldsymbol{a}_1 = T\boldsymbol{u}$ avec $\boldsymbol{u} = (u,v,w)$ transforme les 3 facteurs de la fonction de Plancherel Rotach : $\gamma(\boldsymbol{a}_1)/n = \boldsymbol{s}G(\boldsymbol{a}_1) - Log a_1 bc = L(\boldsymbol{a}_1) + Q(\boldsymbol{a}_1) - Log a_1 bc$ en :

- $Q(\boldsymbol{a}_1)$ en $-\mu v^2 + \mu w^2$
- $L(\boldsymbol{a}_1)$ en $LT\boldsymbol{u} = \mu u + \varpi v/\sqrt{2} + \varpi w/\sqrt{2}$.
- $Log a_1 bc = \log uvw$ car $a_1 bc$ représente le volume du parallélépipède rectangle de côtés $a_1 bc$ que la transformation orthogonale $T$ laisse invariant.

Dans la base $\boldsymbol{u}$, $\gamma(\boldsymbol{u})/n = \gamma_1(u) + \gamma_2(v) + \gamma_3(w)$ est la somme de trois fonctions indépendantes où $\gamma_1(u) = \varpi u/\sqrt{2} - \mu u^2 - \log u$, $\gamma_2(v) = \mu v - \log v$ et $\gamma_3(w) = \varpi w/\sqrt{2} + \mu w^2 - \log w$

D'autre part, comme : $\partial\gamma(\boldsymbol{u})/\partial\boldsymbol{u} = \partial\gamma(\boldsymbol{a}_1)/\partial\boldsymbol{a}_1 \cdot \partial\boldsymbol{a}_1/\partial\boldsymbol{u} = \partial\gamma(\boldsymbol{a}_1)/\partial\boldsymbol{a}_1 \cdot T$, les points critiques de $\gamma(\boldsymbol{a}_1)$ sont ceux de $\gamma(\boldsymbol{u})$ et inversement. Il est donc équivalent de calculer l'écart résolvant $e_f^n(\boldsymbol{a}_1, \boldsymbol{y})\big|_{\boldsymbol{a}_1 = 0} = \partial^n(e^{\boldsymbol{y}\boldsymbol{a}_1} - e^{\boldsymbol{y}f(\boldsymbol{a}_1,\delta)})/\partial \boldsymbol{a}_1^n\big|_{\boldsymbol{a}_1 = 0}$ dans la base $\boldsymbol{a}_1$ ou dans la base $\boldsymbol{u}$.

$sf(a,\delta)$ devient dans la base $\boldsymbol{u}$ : $\gamma(\boldsymbol{u})/n$ dont les dérivées n-ièmes respectives sont :

$\partial^n \exp(\varpi u/\sqrt{2} - \mu u^2)/\partial u^n\big|_{v=0} = H_n(\varpi/2\sqrt{\mu})(-\sqrt{2\mu})^n$, où $H_n$ est le polynôme d'Hermite

$\partial^n \exp(\mu v)/\partial v^n\big|_{u=0} = \mu^n$,

$\partial^n \exp(\varpi w / \sqrt{2} + \mu w^2) / \partial w^n \big|_{w=0} = H_n(\varpi / 2i\sqrt{\mu})(-i\sqrt{2\mu})^n$

$H_n(\varpi / 2i\sqrt{\mu})(i)^n$ est toujours positif sauf si $\varpi = 0$. De même, $\mu^n$ est toujours positif. Les zéros réels de $\partial^n e^{yf(a_1,\delta)} / \partial a_1^n \big|_{a_1=0}$ seront ceux du polynôme d'Hermite $H_n(\varpi / 2\sqrt{\mu})$. Il faut alors établir sous quelles conditions $\partial^n e^{nsf(a_1,\delta)} / \partial a_1^n \big|_{a_1=0}$ l'emporte sur $\partial^n (e^{nsa_1}) / \partial a_1^n \big|_{a_1=0}$ dans la base $\boldsymbol{u}$, autrement dit, pour que le point fixe $\boldsymbol{0}$ ne soit pas attractif. Comme $y a_1 / n = sTu = \mu u + (r/\sqrt{2})(v+w)$, les dérivées 3n-ième de $\exp sTu$ rapport à $\boldsymbol{u}$ sont $\mu^n (r/\sqrt{2})^{2n}$. L'écart résolvant s'écrit enfin : $\mu^n \left[ (r^2/2)^n - H_n(\varpi / 2\sqrt{\mu}) H_n(\varpi / 2i\sqrt{\mu})(i2\mu)^n \right]$

*- Distribution invariante au point fixe $\boldsymbol{0}$*

Le terme de plus haut degré du second membre de l'écart doit être supérieur à $(r^2/2)^n$ pour imposer la solution probabiliste. La condition $|\varpi| > r$ suffit pour que le second membre de l'écart l'emporte et que $\varpi / (2\sqrt[4]{s^2 + t^2})$ soit zéro de $H_n$. On a une famille d'ovales aléatoires de loi bêta $\beta(1/2, 1/2)$,

*- Calcul pour les autres points fixes et distributions dominantes*

Cherchons alors les distributions autour des 2 autres points fixes. Pour passer du point fixe $\boldsymbol{0}$ au point fixe $\boldsymbol{\alpha}_+$ (resp. $\boldsymbol{\alpha}_-$) il suffit de changer $\boldsymbol{a}_1 = (a_1, b, c)$ en $\boldsymbol{a}_1 + \boldsymbol{\alpha}_+ = (a_1, b+\alpha, c+\alpha^2/\beta)$ (resp. $\boldsymbol{a}_1 + \boldsymbol{\alpha}_- = (a_1, b-\alpha, c+\alpha^2/\beta)$) dans l'itération différentielle : $(a + \delta\sigma(b-a)$ ; $b + \delta(\rho a - b - ac)$ ; $c + \delta(-\beta c + ab))$. La proposition permet de passer de $sG(\boldsymbol{a}_1)$ à $sG_+(\boldsymbol{a}_1)$ calculées au point fixe $\boldsymbol{0}$ $sG_+(\boldsymbol{a}_1) = sG(\boldsymbol{a}_1) - (s(\rho-c)+tb)\alpha$ et de $sG(\boldsymbol{a}_1)$ à $sG_-(\boldsymbol{a}_1) = sG(\boldsymbol{a}_1) + (s(\rho-c)+tb)\alpha$.

Donc les deux distributions centrées sur $\boldsymbol{\alpha}_+$ et $\boldsymbol{\alpha}_-$ dominent celle centrée sur $\boldsymbol{0}$. $(s(\rho-c)+tb) = 0$ est voie de communication entre les deux.

*- Interprétation géométrique*

La distribution des zéros du polynôme d'Hermite $H_n(\varpi / 2\sqrt{\mu})$ donne la distribution des ovales aléatoires en $\boldsymbol{0}$. Posons la variable aléatoire $(r + s\rho) / (2\sqrt[4]{s^2 + t^2}) = \chi$. Le plan $c = r + s\rho$ coupe le cylindre $(c/2\chi)^4 = s^2 + t^2$ selon un rayon aléatoire $(c/2\chi)^2$. Pour les distributions centrées sur $\boldsymbol{\alpha}_+$ ou $\boldsymbol{\alpha}_-$, il faut remplacer $\varpi$ par $\varpi_+ = r + s\rho + s\alpha^2/\beta + t\alpha$ ou par $\varpi_- = r + s\rho + s\alpha^2/\beta - t\alpha$ pour obtenir la géométrie asymptotique.

Les premières bases d'une étude mathématique systématique sont jetées bien que cela ne dévalue en rien les approches antérieures. D'où, sous les réserves énoncées :

**Proposition 2**

*Une solution asymptotique de l'équation de Lorenz se compose de 2 familles d'ovales aléatoires $(r + s\rho) / (2\sqrt[4]{s^2 + t^2}) = \chi$ de loi bêta $\beta(1/2, 1/2)$, centrées sur les 2 points fixes $\boldsymbol{\alpha}_+$ et $\boldsymbol{\alpha}_-$, reliées par la voie de communication $(s(\rho-c)+tb) = 0$ et des leurs itérés.*

## 6 - L'équation de Navier Stokes et ses ovales aléatoires

Cet exemple sera traité de façon moins détaillée que celui de Lorenz car on retrouve les mêmes raisonnements avec peu de modifications sinon que $\delta$ n'est plus unidimensionnel.

Fefferman[6], rédacteur du problème de la fondation Clay pour l'équation de Navier Stokes, semble manifestement fort mal à l'aise pour formaliser mathématiquement la question : prudent, il s'en remet aux conceptions du XIX ième siècle où l'on développe les équations différentielles en séries de Fourrier en vue de compactifier la solution. En plus, la présentation est singulièrement prisonnière de considérations physiques, introduisant force paramètres pour rester au plus près de la réalité physique. De crainte qu'en oubliant un, on passe à côté d'une difficulté inaperçue du problème ?

Il faut reconnaître que les mathématiciens sont échaudés par les difficultés que l'on rencontre pour passer de la dimension 2 à 3 (Cela résulte en partie des méthodes utilisées et en partie de la validité du théorème de Bindixon-Poincaré).

Ainsi, les équations de Navier Stokes concernent la vitesse $u$, la pression $p$ comme inconnues d'un fluide où les variables sont la position $x \in \mathbb{R}^n$ et le temps $\tau \in \mathbb{R}^+$ ; on y ajoute la divergence nulle.

On ne s'embarrasse pas ici de telles considérations. On se cale sur nos hypothèses pour déterminer la mesure invariante de Perron Frobenius dans un compact. Il faut d'abord traduire les équations $[6]$ en termes d'itération. Cette opération est alourdie par la multiplicité des équations. On va les mettre sous la forme $\partial a / \partial t = F(a)$ où les inconnues sont notées par $a$, variables sont notées par $t$.

**Traduction de $[6]$ sous la forme $\partial a / \partial t = F(a)$**

- Les inconnues du problème sont le vecteur vitesse $u$ de $\mathbb{R}^n$ de dimension $n$ et la pression $p$.
- Les variables initiales sont la position $x \in \mathbb{R}^n$ et le temps $\tau \in \mathbb{R}^+$.

Au départ, on a autant d'équations que d'inconnues avec la divergence nulle, et l'on va ensuite introduire autant d'équations que d'inconnues intermédiaires :

- La matrice $\partial u / \partial x = b$ introduit les inconnues intermédiaires $b \in \mathbb{R}^n \times \mathbb{R}^n$ avec un nombre égal d'équations. D'où : $\sum_{j=1}^{j=n} u_j \partial u_i / \partial x_j = ab_i$. Avec la matrice $U_1$, composée de zéros ou de 1, on a :

- $\text{div } u = \sum_{j=1}^{j=n} b_{jj} = U_1 b = 0$ ( qui exprime une inconnue intermédiaire en fonction des autres).

- Le champ des forces extérieures $f_i(x, \tau) = f_i(c)$ sera ici polynomial (ou approximé par des polynômes) vérifiant nos hypothèses.

On regroupe sous la dénomination $d$ toutes ces inconnues intermédiaires linéaires aussi variées que $\partial u / \partial x = b$, $\Delta u_i = \sum_{j=1}^{j=n} \partial^2 u_i / \partial x_j^2 = \sum_{j=1}^{j=n} \partial b_{ij} / \partial x_j$ , $\partial p / \partial x_j = d_j$ ou $\partial c / \partial t = 1$ qui se mettent sous la forme $\partial d / \partial t = B$ où $d$ et $B$ sont linéaires en $(a,b,c)$

-Les équations de Navier Stokes $[6]$ s'écrivent en finalement dans la terminologie des itérations :

$\partial a / \partial \tau = -ab + Ad + f(c)$ où $A$ est une matrice à coefficients constants.

$\partial d / \partial t = B$

**Forme de la solution**

On a mis sous la forme $\partial a / \partial t = F(a)$ avec ici $a = (a,b,c,d)$ et $t = (x, \tau)$. Elle est somme toute assez simple puisque quadratique en $a$ et $b$, $A$ est à coefficients constants et $B$ linéaire.

Ainsi, les inconnues sont reliées par des relations linéaires sauf dans les équations principales où ne figurent en définitive que les polynômes $ab$ et $f(c)$. D'où, l'itération $f(a, \delta) = a + \delta F(a)$ où l'on a des linéarités. Supposons qu'elle applique un compact dans lui-même pour tout $\delta > 0$ petit. Formons le $y = (x, y)$ et $P = yF(a,b,c,d) : P = x[-ab + Ad + f(c)] + zc + td + ux$ et $\gamma(a) = P - n \log a$

La hessienne de $P$ est alors , si $^t x$ est le vecteur ligne de $x$ :
$Q = \begin{pmatrix} 0 & ^t x & 0 \\ x & 0 & 0 \\ 0 & 0 & x \partial^2 f / \partial t^2 \end{pmatrix}$

$x \partial^2 f / \partial t^2$ étant calculé en 0. Elle reste constante si $x \partial^2 f / \partial t^2 = 0$. L'équation caractéristique s'écrit : $\mu^k (\mu^2 - |x|^2) |x \partial^2 f / \partial c^2 - \mu I| = 0$ où $k$ est tel que le degré total du polynôme corresponde à la dimension de l'itération. La hessienne est dégénérée et a une valeur propre positive $\mu = |x|$) et une $\mu = -|x|$ ; le mouvement se situe dans espaces propres comme dans le cas de Lorenz. En dimension 3, on a une matrice orthogonale comme T. On a trois groupes de solutions avec des ovales tributaires de la loi bêta $\beta(1/2, 1/2)$, mais aussi des aléas indépendants dus aux forces extérieures. Comme on le voit, la résolution de ces questions implique la maîtrise de bien des paramètres.

En résumé, l'analyse de cette équation diffère fort peu de celle de Lorenz : Seuls sont introduits des points fixes plus nombreux et des aléas éventuels liés à $f$, sinon elle présente le même type d'itération quadratique dégénérée et des linéarités. Elle en diffère toutefois par la présence possible de failles qui peuvent en perturber les comportements. Pour cela, il faut déterminer les valeurs propres de la partie linéaire et les failles qui peuvent en résulter, ce qui ne semble pas très agréable au vu du nombre de matrices binaires introduites. Cela relève davantage du calcul.

### III- Réflexions sur le seuil de récurrence de Poincaré ou l'horizon de Lyapounov ?

**Définition (rappel)**
*Soit $A$ un sous-ensemble mesurable. Un point $a$ est dit **récurrent** par rapport à A si et seulement si pour tout entier $p$, il existe un entier $k \geq p$ pour lequel : $f^k(a) \in A$.*
Le **théorème de Poincaré** affirme alors :
*Si $f(a)$ est bijective et préserve la mesure $\mu$ (pour tout $A \subset \sup p(\mu)$ : $\mu(A) = \mu(f^{-1}(A)) > 0$), alors pour tout $A \subset \sup p(\mu)$, les points $a \in A$ sont presque tous récurrents.*

pour tout $\delta_0 > \delta > 0$ suffisamment petit, $f(a,\delta)$ est inversible assure la récurrence s'il existe une mesure invariante. Si $f(a,\delta)$ applique un domaine compact $C$ de $R^d$ dans lui-même, cette mesure existe et les points du support de la distribution invariante de l'opérateur de Perron-Frobenius sont presque (tous) récurrents.

Supposons que $F(a)$ vérifie nos hypothèses et soit partiellement linéaire de sorte que l'itération différentielle soit équivalente à $G(a) = (a + \tau A(a)b, b)$. Les points critiques $a_c$ et leur répartition sont images réciproques de la répartition uniforme de l'hypercube unité $K = (0,1)^p$ compactifié en ajoutant le point 0. On partage $K$ en cubes $\kappa_\varepsilon$ de côtés $\varepsilon = 1/m$ entourant les points $(1/2, 1+1/2, ..., m-1/2)$.

Prenons alors un pas $\delta$ suffisamment petit pour que $f(a,\delta)$ soit inversible dans le compact. On se place dans la situation où l'on a une mesure invariante. On définit naturellement les variables $t_\delta$ de sorte que les points $a(t_\delta)$ de $f(a,\delta)$ tombent au moins une fois dans chaque ensemble $\alpha$ pourvu que $\delta$ soit très petit : le minimum $t_\delta$ s'interprète comme variable de première visite de chaque cube $\kappa_\varepsilon$. Pour tout $t \ll t_\delta$, l'orbite sera déterministe pour tous ces $t$ en appliquant les théorèmes classiques puisque ses points seront atteints par l'itération différentielle.
Si l'on prend au hasard un sorte que l'ergodisme permette de faire comme si l'on tirait des points au hasard dans $K$, alors le comportement sera aléatoire en appliquant Perron Frobenius. Seuls seront définis la densité de probabilité et son support : la position n'a plus de sens.
On voit que cette conception, proche de la mécanique quantique précise la stabilité au sens de Lyapounov. (voir Ghys $[8]$). Le théorème de Bendixon-Poincaré interdit ce cas en dimension 1 ou 2.

**Remarques concernant les mathématiques**
- Certaines solutions asymptotiques d'une équation différentielle dans conditions particulières peuvent donc être aléatoires. Le modèle probabiliste apporte alors les informations que l'approche déterministe est incapable de maîtriser.
- Chose très curieuse, la distribution asymptotique pour $F(a)$ donné est la même, qu'il s'agisse de EDO ou de EDP, mais les conditions de convergence sont très différentes. En particulier pour les EDP, une petite variation sur la direction $\tau$ peut faire basculer le comportement de l'itération différentielle soit vers un point fixe (ou cycle) soit vers un aléa. Reste à étudier le détail des cycles.
- Beaucoup d'idées sur les chaos doivent être corrigées : la définition du chaos dans la littérature, comme situation intermédiaire entre déterminisme et probabilité, n'a pas de sens. L'idée que les multiples erreurs d'approximation rendent l'évolution imprévisible cache les vraies raisons mathématiques. Enfin, la notion de sensibilité aux conditions initiales doit être complètement revue.
- En revanche, la dégénérescence de la hessienne prend une importance considérable. Cette dégénérescence reste encore un problème délicat à traiter. Par exemple, si l'on a l'équation $da/dt = F(a,t)$ au lieu de $da/dt = F(a)$, il suffit d'ajouter une coordonnée $dt/dt' = 1$ et de prendre $t'$ comme variable, ce qui introduit une linéarité : la dégénérescence de la hessienne modifie la structure des résultats et explique la complexité de ces phénomènes non autonomes.

**Remarques concernant la physique**

- Les conditions aux limites sont alors essentielles pour assurer la compacité et de la récurrence de l'itération en physique pour remplir la condition : tout voisinage d'un point du domaine devra être « visité » un nombre infini de fois par l'itération lorsque $t \to \infty$. Ainsi, dans la plupart des problèmes de la physique, on s'arrange pour que des conditions aux limites concernant $t = (q,t')$ mettent une borne aux variations de $a$ pour les réduire à un domaine compact. Dans certains problèmes, le domaine de $q$ a un contour borné fini. Le travail de Delabaère et Howls [5] permet d'introduire ces conditions sans trop de difficulté dans notre présentation. Dans d'autres problèmes, la périodicité spatiale (par exemple l'équation de Navier-Stokes [6]) de la solution $a(q,t') = a(q+e,t')$ est la plus sure contrainte de récurrence, seul le temps pouvant aller de son libre cours.
- L'hypothèse de compacité concerne la plupart des applications de la physique à distance finie.
- Enfin, l'idée que de multiples erreurs d'approximation rendent l'évolution imprévisible, idée que l'on a combattue précédemment en mathématique, a tout son sens en physique.

**5- Exemple** : **Les équations d'Hamilton dans $R^d$**

- En gardant les notations de la mécanique classique [1] :

$dp/dt = -\partial H(p,q)/\partial q$ ; $dq/dt = \partial H(p,q)/\partial p$.

Ces équations réalisent une application de $R^{2d}$ dans $R^{2d}$ quand l'hamiltonien $H$ est autonome. Associons à $F = (-\partial H(p,q)/\partial q, \partial H(p,q)/\partial p)$ l'itération $f(\boldsymbol{a},\delta) = \boldsymbol{a} + \delta F(\boldsymbol{a})$, $\boldsymbol{a} = (p,q)$ et supposons qu'elle applique un compact dans lui-même pour tout $\delta > 0$ petit.

Il faut rechercher les points fixes (de Lagrange) $\partial H/\partial q = 0$ et $\partial H/\partial p = 0$, puis, calculer les valeurs propres de la partie linéaire pour ces points fixes, en posant $I$ la matrice identité de $R^d$ et en annulant le déterminant formé des 4 blocs $d \times d$ :

$$J(\lambda) = \begin{vmatrix} -\partial^2 H/\partial q^2 - \lambda I & -\partial^2 H/\partial p \partial q \\ \partial^2 H/\partial p \partial q & \partial^2 H/\partial p^2 - \lambda I \end{vmatrix} = 0$$

Soit $\gamma(p,q) = -x\partial H(p,q)/\partial q + y\partial H(p,q)/\partial p - n\log p - n\log q$ la fonction de Plancherel-Rotach ; la densité invariante devra être recherchée parmi points critiques imaginaires :

$n/q = -x\partial^2 H/\partial q^2 + y\partial^2 H/\partial p \partial q$

$n/p = -x\partial^2 H/\partial q \partial p + y\partial^2 H/\partial p^2$

définissant en général la distribution intrinsèque en cas de dégénérescence de la hessienne $\partial^2 \gamma/\partial q \partial p$ sauf pour les cycles où la distribution est uniforme avec $\int_0^T \partial H/\partial q \, dt = 0$ et $\int_0^T \partial H/\partial p \, dt = 0$.

En résumé, l'analyse de cette équation diffère fort peu de celle de Lorenz : Seuls sont introduits les points fixes et les aléas liés à $f$. L'analyse en diffère toutefois par la présence possible de failles qui peuvent venir perturber les comportements. Il reste donc à déterminer les valeurs propres de la partie linéaire et les failles qui peuvent en résulter.

## III- Équations intégrales avec noyau d'Écalle-Legendre [7]

Certaines équations intégrales peuvent être traitées à partir du lemme 1 du premier article [3].

La transformée de Laplace $\phi_f(y)$ de $p_f$ associée à l'itération $f$ et à la densité $p$ arbitraires est :

$$\phi_f(y) = \int p(s)e^{yf(s)}ds$$

Or, $p$ est l'inverse de $\phi$ : $p(s) = \frac{1}{(2\pi i)^d}\int_{c-i\infty}^{c+i\infty} e^{-xs}\phi(x)dx$. Si $N(y,x)$ est le noyau d'Écalle défini par $N(y,x) = \int e^{yf(s)-xs}ds$, l'équation résolvante $\phi_f(y) = \phi(y)$ se traduit par $\phi(y) = \frac{1}{(2\pi i)^d}\int_{c-i\infty}^{c+i\infty} N(y,x)\phi(x)dx$ et permet de trouver des solutions probabilistes de l'équation intégrale $\phi(y) = \frac{1}{(2\pi i)^d}\int_{c-i\infty}^{c+i\infty} N(y,x)\phi(x)dx$. On peut trouver une approximation normale de ce noyau.

**IV- Les équations avec décalage $a(\lambda t) = f \circ a(t)$ ou $f \circ a(t) = \lambda a(t)$**

Ces équations déterminent les courbes semi-invariantes des itérations de $f(a)$ où $\lambda$ sont les valeurs propres de la partie linéaire de $f$. Elles peuvent présenter des structures des plus curieuses, en particulier celle relative à l'itération de Henon que l'on a particulièrement étudiée sur ordinateur.

**Bibliographie**